\documentclass{amsart}

\newtheorem{theorem}{Theorem}[section]
\newtheorem{lemma}[theorem]{Lemma}
\newtheorem{proposition}[theorem]{Proposition}

\theoremstyle{definition}
\newtheorem{definition}[theorem]{Definition}
\newtheorem{example}[theorem]{Example}

\theoremstyle{remark}
\newtheorem{remark}[theorem]{Remark}

\makeatletter
\newcommand{\cp}{\mathop{\operator@font cp}}
\newcommand{\range}{\mathop{\operator@font range}}
\newcommand{\rank}{\mathop{\operator@font rank}}
\newcommand{\dom}{\mathop{\operator@font dom}}
\newcommand{\Real}{\mathop{\operator@font Re}}
\newcommand{\cont}{\mathop{\operator@font Cont_w}}
\makeatother

\begin{document}
\title{Compact Operators in TRO's}
\author{G. Andreolas}
\address{Department of Mathematics, University of the Aegean, 832\,00
Karlovassi, Samos, Greece}
\email{gandreolas@aegean.gr}
\date{}


\begin{abstract}
 We give a geometric characterization of the elements of a TRO that can be represented as
compact operators by
a faithful representation of the TRO.
\end{abstract}

\maketitle

\section{Introduction}

A \textit{ternary ring of operators} (or simply, TRO) between Hilbert spaces $H_2$ and
$H_1$ is a norm closed subspace $\mathcal{V}$ of $\mathcal{B}(H_2,H_1)$ which is closed
under the triple
product

$$\mathcal{V}\times \mathcal{V}\times \mathcal{V} \ni(x,y,z) \mapsto xy^*z\in
\mathcal{V}.$$
A TRO $\mathcal{V}\subseteq \mathcal{B}(H_2,H_1)$ is called a w*-TRO if it is w*-closed
(equivalently, weak
operator closed, or strong operator closed) in
$\mathcal{B}(H_2,H_1)$. TRO's were first introduced by Hestenes \cite{hes} and since then
they have
been studied by many authors.
In general, a TRO $\mathcal{V}$ can be identified with the off-diagonal corner (at the
(1,2)
position) of its \textit{linking} C*-algebra

$$A(\mathcal{V})= 
 \left(
   \begin{array}{ll}
    \mathcal{C} & \mathcal{V}\\
    \mathcal{V}^* & \mathcal{D}
   \end{array} \right)
\subseteq \mathcal{B}(H_1\oplus H_2),$$
where $\mathcal{C}$ and $\mathcal{D}$ are the C*-algebras generated by
$\mathcal{V}\mathcal{V}^*$ and $\mathcal{V}^*\mathcal{V}$ respectively.

If $\mathcal{S}$ is a nonempty subset of the unit ball of a normed space $\mathcal{X}$,
then
the \textit{contractive perturbations} of $\mathcal{S}$ are defined as
$$\cp(\mathcal{S})=\left\{x\in \mathcal{X}\ |\ \|x\pm s\|\leq 1\ \forall s\in
\mathcal{S}\right\}.$$
It is clear that $\mathcal{S}_1\subseteq \mathcal{S}_2$ implies
$\cp(\mathcal{S}_1)\supseteq \cp(\mathcal{S}_2)$. Also, an element $x$
of the unit ball of $\mathcal{X}$ is an extreme point if
and only if $\cp(\{x\})=\{0\}$. We shall write $\cp(x)$ instead of $\cp(\{x\})$. 

One may define contractive perturbations of higher-order by using the recursive
formula $\cp^{n+1}(\mathcal{S})=\cp\left(\cp^{n}(\mathcal{S})\right),$ $n\in\mathbb N$. It
is clear that
$\cp(\mathcal{S})$ is a norm-closed convex subset of the closed unit ball of
$\mathcal{X}$. One can also verify
that
$\mathcal{S}\subseteq\cp^{2}(\mathcal{S})$; from this it follows that
$\cp^{3}(\mathcal{S})=\cp(\mathcal{S})$. The second
contractive perturbations were introduced in \cite{1996}. In \cite{1996} it is
proved that the set of the second contractive perturbations of an element $a$ of
a C*-algebra $\mathcal{A}$ is compact in the norm topology
if and only if there exists a faithful representation $\phi$ of
$\mathcal{A}$ such that $\phi(a)$ is a compact operator. Further study was
conducted in \cite{2008}, \cite{1997}, \cite{2005} and \cite{2004}. We shall see that this
characterization is not valid for the elements of a TRO.

In this work we characterize the elements of a TRO that are represented as compact
operators by a faithful representation of the TRO, in terms of the size of
their contractive perturbations. We show that
there exists a
faithful representation $\phi$ of the TRO $\mathcal{V}$ that maps an element $a$ of the
unit ball of $\mathcal{V}$ to a compact operator if and only if the
set of its second contractive perturbations is weakly compact. It follows from
\cite{1996} and our result that for an element
$a$ of a C*-algebra $\mathcal{A}$ the set $\cp^2(a)$ is compact if and only if it is
weakly compact or, equivalently, there exists a faithful representation $\pi$ of
$\mathcal{A}$ such that $\pi(a)$ is a compact operator. 

Ylinen proved in \cite{1972} and
\cite{1975} that for an element $a$ of a C*-algebra $\mathcal{A}$ the operator
$x\rightarrow axa$ on $\mathcal{A}$ is compact if and only if it is weakly compact or,
equivalently, there exists a faithful representation $\pi$ of $\mathcal{A}$ such that
$\pi(a)$ is a compact operator. We obtain an analogous result for the operator
$x\rightarrow ax^*a$ on a TRO.

\vspace{1em}
\noindent \textbf{Notation.} Throughout, we adopt the following notation: $H_1$
and $H_2$ are Hilbert spaces, $\mathcal{B}(H_2,H_1)$ the space of all bounded linear
operators $H_2\rightarrow\! H_1$ and $\mathcal{K}(H_2,H_1)$ the space of all compact
operators $H_2\rightarrow H_1$. In
particular $\mathcal{B}(H_1)=\mathcal{B}(H_1,H_1)$ and
$\mathcal{K}(H_1,H_1)=\mathcal{K}(H_1)$. $\mathcal{V}$ is a TRO that is a
subspace of $\mathcal{B}(H_2,H_1)$. Let $\mathcal{X}$ be a Banach space,
$\mathcal{Y}\subseteq \mathcal{X}$ a subspace and
$a\in \mathcal{Y}$. Then by $\cp^n_{\mathcal{Y}}(a)$ we denote the set of the n-th
contractive
perturbations of $a$ computed with respect to $\mathcal{Y}$. If $r$ is a positive
number, then by
$\mathcal{X}_r$ we denote the closed ball of center 0 and radius $r$. Let $x,y$ be
elements of a
Hilbert space $H$. We denote by $x\otimes y$ the rank one operator on $H$ defined by
$$(x\otimes y)(z)=\langle z,x\rangle y.$$

\section{Preliminaries}

Let $\mathcal{V}$ and $\mathcal{W}$ be two TRO's. A linear map
$\phi:\mathcal{V}\rightarrow
\mathcal{W}$ is
called a
\textit{TRO-homomorphism} if it preserves the ternary product

$$\phi(xy^*z)=\phi(x)\phi(y)^*\phi(z)$$
for all $x,y,z\in \mathcal{V}$. If, in addition $\phi$ is an injection from $\mathcal{V}$
onto $\mathcal{W}$, we call
$\phi$ a \textit{TRO-isomorphism} from $\mathcal{V}$ onto $\mathcal{W}$. A
TRO-homomorphism
$\phi$ from a TRO $\mathcal{V}$ into the set of all bounded operators from one Hilbert space
to
another, is called a \textit{representation} of $\mathcal{V}$. We will say that a
representation
$\phi:\mathcal{V}\rightarrow \mathcal{B}(H_2,H_1)$ is a \textit{faithful} representation
of $\mathcal{V}$ if
$\phi$ is injective. It was shown in \cite[Proposition 3.4]{1981} that every faithful
TRO-representation is an isometry.

\begin{proposition} \label{incl1}
 Let $H_1$ and $H_2$ be Hilbert spaces, $\mathcal{V}\subseteq \mathcal{B}(H_2,H_1)$ a TRO
and
$A(\mathcal{V})$ its linking
algebra.  If $a$ is in the unit ball of $\mathcal{V}$, then 
$\cp_{A(\mathcal{V})}^2(a)\subseteq \cp_\mathcal{V}^2(a)$.
\end{proposition}
\begin{proof}
First we note that if $\mathcal{E}$ is a Banach space and $b\in \mathcal{E}$ has the
property $\|x+b\|\leq 1$
for all $x\in \mathcal{E}$ with $\|x\|\leq 1$, then $b=0$. Indeed, if the above property
holds for
some $b\neq 0$, then taking $x=b/\|b\|$ we have $\|b/\|b\|+b\|\leq 1$ which implies
$\|b\|=0$. This yields a contradiction.

Let 
$$\left(
   \begin{array}{cc}
    b_1 & b_2\\
    b_3 & b_4
   \end{array} \right)
\in {\cp}_{A(\mathcal{V})}^{2}(a).$$
We can easily see that 
$$\left(
   \begin{array}{cc}
    0 & x\\
    y & 0
   \end{array} \right)
\in {\cp}_{A(\mathcal{V})}(a)$$
for every  $x\in \cp_\mathcal{V}(a)$ and $y\in \mathcal{V}^*$ with $\|y\| \leq 1$. So, it
follows
directly that $b_2\in \cp_{\mathcal{V}}^2(a)$ while from the remark at the beginning of
the proof it
follows that $b_3=0$. We have that $\left(
   \begin{array}{cc}
    b_1 & b_2\\
    y & b_4
   \end{array} \right)$
is a contraction for every $y\in \mathcal{V}^*$ with $\|y\|\leq 1$. Thus, if $\eta\in
H_1$, we have that
\begin{equation} \label{inequality}
 \|b_1\eta\|^2+\|y\eta\|^2\leq \|\eta\|^2
\end{equation}
for all $y\in \mathcal{V}^*$ with $\|y\|\leq1$. Since the strong*-topology of
$\mathcal{V}^*$ is finer than
its strong topology, it follows form \cite[Theorem 3.6 (Kaplansky density
theorem)]{1981} that the
inequality (\ref{inequality}) holds for all $y$ in the closed unit ball of
$\overline{\mathcal{V}^*}^{w^*}$.
Therefore, for
all partial isometries $y\in \overline{\mathcal{V}^*}^{w^*}$, we have $0\leq
b_1^*b_1+y^*y\leq1$.
Denoting by $p_y$ the domain projection $y^*y$ of a partial isometry $y\in
\overline{\mathcal{V}^*}^{w^*}$, it follows that $0\leq b_1^*b_1\leq 1-p_y$.
Multiplying by
$p_y$, we deduce that $p_yb_1^*b_1p_y=0$ or $b_1p_y=0$. Let $\Pi$ be the set of all
partial isometries of $\overline{\mathcal{V}^*}^{w^*}$ and $p$ the orthogonal projection
onto the closed linear span of the subspaces
$\{p_y(H_2)\}_{y\in \Pi}$. Then we have proved that
$b_1p=0$. On the other hand, we can see that
$b_1p^{\perp}=0$, since $b_1$ is in the C*-algebra generated by $\mathcal{V}\mathcal{V}^*$
and
$\overline{\mathcal{V}^*}^{w^*}$ is generated by its partial isometries \cite[Theorem
3.2]{1981}.
Hence, we have
proved that
$b_1=0$. By symmetry, we obtain $b_4=0$. Thus, we
showed that each element of $\cp_{A(\mathcal{V})}^2(a)$ is in $\mathcal{V}$. The fact that
$\cp_{A(\mathcal{V})}^2(a)\subseteq \cp_{\mathcal{V}}^2(a)$ is immediate.
\end{proof}

\begin{remark} \label{strict}
 The containment in the last proposition may be strict. We shall give an example. Let
$H_1$ and $H_2$ be Hilbert spaces with $\dim H_1=\infty$ and $\dim H_2<\infty$ and
$u:H_2\rightarrow
H_1$ an isometry. Let $\mathcal{V}=\mathcal{B}(H_2,H_1)$. Since $u$ is an extreme point of
$\mathcal{V}$ \cite{zettl}, the set $\cp_{\mathcal{V}}^2(u)$
is equal to the unit ball of $\mathcal{V}$. Now, $A(\mathcal{V})=\mathcal{B}(H_1\oplus H_2)$ and it follows from
\cite[Corollary 2.4]{1996} that $\cp_{A(\mathcal{V})}^2(u)$ is compact. Hence, the inlusion
$\cp_{A(\mathcal{V})}^2(u)\subset \cp_\mathcal{V}^2(u)$ is strict. Considering the identity representation of
$\mathcal{V}$ in this example, one can see that the implication $(i)\Rightarrow (ii)$ of \cite[Theorem 2.2]{1996} does
not hold for TRO's.
\end{remark}

\begin{remark} 
It is known that the linking algebra $A(\mathcal{V})$ is
just the C*-envelope $C^*_e(\mathcal{V})$ of the
TRO $\mathcal{V}$. Therefore, the inclusion in Proposition \ref{incl1} in the case of an
operator space $\mathcal O$ would be
$\cp^2_{C^*_e(\mathcal O)}(a)\subseteq \cp^2_{\mathcal O}(a)$. Now, we shall see that this
inclusion does not hold in
operator spaces in general. 

Let $H$ be an infinite dimensional Hilbert space,
$$\mathcal{O}=\left\{\left[
\begin{array}{cc}
    \lambda\, Id & a\\
    0 & \mu\, Id
   \end{array} \right]
\ :\ a\in \mathcal{K}(H), \lambda, \mu\in\mathbb C\right\}.$$
The C*-algebra generated by $\mathcal{O}$ in $\mathcal{B}(H\oplus H)$ is 
$$C^*_{\mathcal{B}(H\oplus H)}(\mathcal{O})=\left\{\left[
\begin{array}{cc}
    \lambda\, Id+a & b\\
    c & \mu\, Id+d
   \end{array} \right]
:\ a, b, c, d\in \mathcal{K}(H), \lambda, \mu\in\mathbb C\right\}.$$
If $\mathcal{I}$ is a proper ideal of $C^*_{\mathcal{B}(H\oplus H)}(\mathcal{O})$, then
$\mathcal{I}$ contains $\mathcal{K}(H\oplus
H)$, the compact operators on $H\oplus H$ and, consequently, the quotient space
$C^*_{\mathcal{B}(H\oplus H)}(\mathcal{O})/\mathcal{I}$ is finite
dimensional. Therefore, $C^*_{e}(\mathcal{O})=C^*_{\mathcal{B}(H\oplus H)}(\mathcal{O})$.
Then if we consider the operator
$$s=\left[
\begin{array}{cc}
    Id & 0\\
    0 & 0
   \end{array} \right],$$
we see that 
$${\cp}^2_{C^*_e}(s)=\left\{\left[
\begin{array}{cc}
    \lambda\, Id+x & 0\\
    0 & 0
   \end{array} \right]:\ 
\lambda\in\mathbb C,\, x\in \mathcal{K}(H),\, \|\lambda\, Id+x\|\leq 1\right\},$$
and that
$${\cp}^2_\mathcal{O}(s)=\left\{\left[
\begin{array}{cc}
    \lambda & 0\\
    0 &  0
   \end{array} \right]:\ 
\lambda\in\mathbb C,\, |\lambda|\leq 1\right\}.$$

\end{remark}

\begin{proposition}\label{cp}
Let $H_1$ and $H_2$ be Hilbert spaces with $\dim H_1=\infty$ and $\dim
H_2=\infty$ and $\mathcal{V}\subseteq \mathcal{B}(H_2,H_1)$ a TRO.
 Let $a=\sum_{i=1}^{\infty}\lambda _i u_i \in \mathcal{V}$ be  a norm one
compact operator, where $\{u_i\}_{i=1}^{\infty}$ are finite rank partial isometries such
that $u_iu_j^*=0$, $u_i^*u_j=0$ for $i\neq j$ and
$\{\lambda_i\}_{i=1}^{\infty}$ is a sequence of positive numbers decreasing to $0$.
Let $e_k=\sum_{i=1}^{k}u_i^*u_i$ and $f_k=\sum_{i=1}^{k}u_iu_i^*$. If $x$ is any
contraction in $\mathcal{V}$, then 
\begin{enumerate}
\item $\|a\pm (1-\lambda_k)f_k^{\bot}xe_k^{\bot}\|\leq 1,$
\item $e_k$ and $f_k$ are in the C*-algebra generated by $\mathcal{V}^*\mathcal{V}$ and
$\mathcal{V}\mathcal{V}^*$ respectively,
\item $f_k^{\bot}xe_k^{\bot}\in \mathcal{V}$.
\end{enumerate}
\end{proposition}
\begin{proof}
(1) Let $y$ be a contraction in $\mathcal{V}$. From \cite{1987} we know that
$$\|a\pm (\textbf{1}-|a^*|)^{1/2}y(\textbf{1}-|a|)^{1/2}\|\leq 1.$$
Simple computations show that 
$$\left\|a\pm \left(\sum_{i=1}^{\infty}(1-\lambda _i)^{1/2}u_i u_i^*
+f^{\bot}\right)y\left(\sum_{j=1}^{\infty}(1-\lambda_j)^{1/2}u_j^*
u_j+e^{\bot}\right)\right\|\leq 1,$$
where $e=[\{u_i^*u_i\}_{i\in \mathbb N}]$ and $f=[\{u_iu_i^*\}_{i\in \mathbb N}]$. Now
setting 
\begin{multline*}y=\left(\sum_{i=k+1}^{\infty}\frac{(1-\lambda_k)^{1/2}}{(1-\lambda_i)^{
1/2}}u_iu_i^*+(1-\lambda_k)^{1/2}f^{\bot}\right)x \\
\left(\sum_{j=k+1}^{\infty}\frac{(1-\lambda_k)^{1/2}}{(1-\lambda_j)^{1/2}}
u_j^*u_j+(1-\lambda_k)^{1/2}e^{\bot}\right),
\end{multline*}
where $x\in \mathcal{V}$ is a contraction, we get the result.

(2) Assume that $\lambda_1=1$. We define a sequence $(a_i)_{i\in \mathbb N}$ in
$\mathcal{V}$, where
$a_1=a$ and $a_n=a_{n-1}a_{n-1}^*a_{n-1}$. Simple computations show that
$\lim_{n\rightarrow\infty}a_n=u_1$ is in $\mathcal{V}$. That means
$a-u_1=\sum_{i=2}^{\infty}\lambda_iu_i$ is in $\mathcal{V}$ and using the same argument
for
$a-u_1$ we deduce $u_2\in \mathcal{V}$ and continuing in the above fashion, we
inductively get $u_n \in \mathcal{V}$ for all $n\in \mathbb N$. Hence, $u_n^* \in
\mathcal{V}^*$ for all $n\in \mathbb N$. It follows that $e_k=\sum_{i=1}^k
u_i^*u_i = (\sum_{l=1}^k u_l^*) (\sum_{m=1}^k u_m)\in \mathcal{V}^*\mathcal{V}$ and
$f_k=\sum_{i=1}^k u_iu_i^* = (\sum_{l=1}^k u_l) (\sum_{m=1}^k u_m^*)\in
\mathcal{V}\mathcal{V}^*$.

(3) Since $x\in \mathcal{V}$, $xe_k\in \mathcal{V}$, $f_kx\in \mathcal{V}$ and $f_kxf_k\in
\mathcal{V}$, it follows
that $f_k^{\bot}xe_k^{\bot}=(\textbf{1}-f_k)x(\textbf{1}-e_k)=x-xe_k-f_kx+f_kxe_k$ is in
$\mathcal{V}$.
\end{proof}

\begin{proposition}\label{incl}
Let $H_1$ and $H_2$ be Hilbert spaces with $\dim H_1=\infty$ and $\dim
H_2=\infty$ and $\mathcal{V}\subseteq \mathcal{B}(H_2,H_1)$ a TRO.
If \;$\mathcal{C}$ is the TRO that consists of all compact operators of $\mathcal{V}$,
then
$\cp_\mathcal{V}^2(a)\subseteq \cp_\mathcal{C}^2(a)$ for all $a\in \mathcal{C}$.
\end{proposition}
\begin{proof}
 Let $a\in \mathcal{C}$. It suffices to show that $\cp_\mathcal{V}^2(a)\subseteq
\mathcal{K}(H_2,H_1)$. We shall show
that if $x\in \mathcal{V}\backslash \mathcal{C}$, then $x\notin \cp_\mathcal{V}^2(a)$.
Since the operator $x$ is not
compact, there exists an $\varepsilon >0$ such that for all finite rank projections $f$,
$e$ on $H_1$, and $H_2$ respectively, the inequality
$$\|f^{\perp}xe^{\perp}\|>\varepsilon$$
holds, where $f^{\perp}=\textbf{1}-f$ and $e^{\perp}=\textbf{1}-e$.
Given that the operator $a$ is compact, there exists a unique sequence of positive numbers $(\lambda_i)_{i\in\mathbb
N}$ desceasing to 0 and a sequence
$\{u_i\}_{i=1}^{\infty}$ of finite rank partial isometries with $u_iu_j^*=0$,
$u_i^*u_j=0$ for $i\neq j$ such that
$$a=\sum_{i=1}^{\infty}\lambda_i u_i.$$ 

Let $e_k=\sum_{i=1}^{k}u_i^*u_i$ and $f_k=\sum_{i=1}^{k}u_iu_i^*$ for all 
$k\in\mathbb N$. From Proposition \ref{cp} we know that if $y$ is any contraction in
$\mathcal{V}$,
then $(1-\lambda_k)f_k^{\perp}ye_k^{\perp}\in \cp_\mathcal{V}(a)$. Thus, it suffices to
find a $k\in
\mathbb N$ and a contraction $y\in \mathcal{V}$ such that $\|x\pm
(1-\lambda_k)f_k^{\perp}ye_k^{\perp}\|>1$.

 We choose $k$ so that $\lambda_{k}<\varepsilon$, set 
$x_{k}=f_{k}^{\perp}xe_{k}^{\perp}$ and $y=x_{k}/\|x_{k}\|$.
The following computations complete the proof
$$\|x+(1-\lambda_{k})f_{k}^{\perp}ye_{k}^{\perp}\|\geq\|x_{k}+(1-\lambda_{k})y
\|=$$
$$\left\|x_{k}+(1-\lambda_{k})\frac{x_{k}}{\|x_{k}\|}\right\|=
\|x_{k}\|\left|1+(1-\lambda_{k})\frac{1}{\|x_{k}\|}\right|=$$
$$\|x_{k}\|+(1-\lambda_{k})>\varepsilon+1-\varepsilon=1.$$
\end{proof}

\section{The Main Results}

We have seen in Remark \ref{strict} that the characterization given in \cite[Theorem 2.2]{1996} does not
hold for TRO's. In this section we shall show that there exists a faithful representation
$\phi$ of the TRO $\mathcal{V}$ that maps an element $a\in
\mathcal{V}_1$ to a compact operator, if and only if the set $\cp_\mathcal{V}^2(a)$ is
weakly compact. This is one of the main results of this work.

Note that if $\pi$ is a faithful representation of a TRO $\mathcal{V}$, we can identify
$\mathcal{V}$ with $\pi(\mathcal{V})$.

\begin{lemma}\label{sa}
Let $a$ be a non-compact selfadjoint operator in $\mathcal{B}(H)_1$. Then there
exists $\varepsilon > 0$ and an infinite dimensional projection $p$ on $H$ such that
$\mathcal{B}(p(H))_{\varepsilon^2/2}\subseteq a\mathcal{B}(H)_{1/2}a$.
\end{lemma}
\begin{proof}
Let us assume that $a$ is a non-compact selfadjoint contraction. We shall denote by
$E$ the unique spectral measure relative to $(\sigma(a),H)$ such that $a=\int
zdE$, where
$z$ is the inclusion map of $\sigma(a)$ in $\mathbb C$. From \cite[Proposition 4.1]{cnw}
there exists an $\varepsilon>0$ such that the projection $p=E(\{z\in\sigma
(a):|z|>\varepsilon\})$ is infinite dimensional. Denote by $a_p$ the operator in
$\mathcal{B}(p(H))$
such that $a_p(h)=ap(h)=pa(h)$ for all $h\in p(H)$. The operator
$a_p$ is invertible. Let us assume that the operator $T$ is in
$p\mathcal{B}(H)_{\varepsilon^2/2}p = \mathcal{B}(p(H))_{\varepsilon^2/2}$. Then
$$\|(a_p)^{-1}T(a_p)^{-1}\|\leq \|(a_p)^{-1}\|^2\|T\|\leq
\frac{1}{\varepsilon^2}\frac{\varepsilon^2}{2}=\frac{1}{2}.$$
Therefore, 
$$T=a_p((a_p)^{-1}T(a_p)^{-1})a_p\in a_p\mathcal{B}(p(H))_{1/2}a_p\subseteq
ap\mathcal{B}(H)_{1/2}pa.$$
So,
$$\mathcal{B}(p(H))_{\varepsilon^2/2}=p\mathcal{B}(H)_{\varepsilon^2/2}p\subseteq
ap\mathcal{B}(H)_{1/2}pa\subseteq a\mathcal{B}(H)_{1/2}a.$$
\end{proof}

\begin{proposition} \label{bofh}
 Let $a$ be a contractive operator on a Hilbert space $H$. Then the operator $a$ is
compact if and only if the set $\cp_{\mathcal{B}(H)}^2(a)$ is weakly compact.
\end{proposition}
\begin{proof}
 The forward implication is trivial from \cite[Corollary 2.4]{1996}.

Conversely, suppose that the operator $a$ is non-compact. The polar decomposition of
$a$ has the following form
$$a=v|a|,$$
where $v$ is a partial isometry, such that $v^*v|a|=|a|$ and
$\dom(v)=\overline{|a|(H)}$.
From Lemma \ref{sa} we know that there exists $\varepsilon>0$ and an
infinite dimensional projection $p$ such that 
$$vp\mathcal{B}(H)_{\varepsilon^2/2}p\subseteq v|a|\mathcal{B}(H)_{1/2}|a|.$$
Therefore, the following inclusions hold
$$\mathcal{B}(p(H),vp(H))_{\varepsilon^2/2}= vp\mathcal{B}(H)_{\varepsilon^2/2}p \subseteq
v|a|\mathcal{B}(H)_{1/2}|a|=$$
$$v|a|\mathcal{B}(H)_{1/2}v^*v|a|\subseteq
v|a|\mathcal{B}(H)_{1/2}v|a|=a\mathcal{B}(H)_{1/2}a\subseteq{\cp}^2_{\mathcal{B}(H)}(a).$$
The last inclusion follows from \cite[Proposition 1.2]{1996}.
Since $v$ is a non-compact partial isometry, $\mathcal{B}(p(H),vp(H))_{\varepsilon^2/2}$
is not
weakly compact \cite[Chapter V, Theorem 4.2]{cnw}. The proof is complete.
\end{proof}

\begin{theorem} \label{cstarweakly}
 Let $\mathcal{A}$ be a C*-algebra and $a\in \mathcal{A}_1$. Then there exists a faithful
representation
$\phi$ of $\mathcal{A}$ such that $\phi(a)$ is a compact operator if and only if
$\cp^2(a)$
is a weakly compact set.
\end{theorem}
\begin{proof}
 The forward implication follows from \cite[Theorem 2.2]{1996}.

Conversely assume that $\phi(a)$ is a non-compact operator for all faithful
representations $\phi$ of $\mathcal{A}$. Let $\{(\phi_i,H_i)\}$ be a
maximal family of pairwise inequivalent irreducible representations of $\mathcal{A}$ and
let $\phi$
be the reduced atomic representation $(\phi,\sum_{i\in I}\oplus H_i)$. Since all $\phi_i$
are irreducible representations, the SOT-closure of $\phi(\mathcal{A})$ equals $\sum_{i\in
I}\oplus
\mathcal{B}(H_i)$. Kaplansky's Density Theorem shows that $\phi(\mathcal{A}_1)$ is
SOT-dense in $(\sum_{i\in
I}\oplus \mathcal{B}(H_i))_1$ and so $\phi(a)\phi(\mathcal{A}_{1/2})\phi(a)$ is SOT-dense
in
$$\sum_{i\in I}\oplus\phi_i(a)\mathcal{B}(H_i)_{1/2}\phi_i(a).$$
However, \cite[Proposition 1.2]{1996} shows that $\phi(a)\phi(\mathcal{A}_{1/2})\phi(a)$
is
contained in the set $\cp^2(\phi(a))$, which is a SOT-closed set. Thus
$$\left(\sum_{i\in I}\oplus \phi_i(a)\mathcal{B}(H_i)_{1/2}\phi_i(a)\right)\subseteq
{\cp}^2(\phi(a)).$$
The operator $\phi(a)$ is not compact,
since the reduced atomic representation is faithful. Thus, there are two cases.

Assume first that there exists an $i_o\in I$ such that $\phi_{i_o}(a)$ is a
non-compact operator on $H_{i_o}$. Therefore, from the proof of Proposition \ref{bofh}
there exists an infinite dimensional projection $p\in \mathcal{B}(H_{i_o})$, a non-compact
partial
isometry $v$ and an $\varepsilon>0$ such that
$\mathcal{B}(p(H_{i_o}),vp(H_{i_o}))_{\varepsilon^2/2}\subseteq
\phi_{i_o}(a)\mathcal{B}(H_{i_o})_{1/2}\phi_{i_o}(a)$. It follows that
$\mathcal{B}(p(H_{i_o}),vp(H_{i_o}))_{\varepsilon^2/2}\oplus\sum_{i\in I-\{i_o\}}\oplus
\phi_i(a)\mathcal{B}(H_i)_{1/2}\phi_i(a)\subseteq \cp^2(a)$. Therefore the set
$\cp^2(\phi(a))$ is
not weakly compact since $\mathcal{B}(p(H_{i_o}),vp(H_{i_o}))_{\varepsilon^2/2}$ is not a
weakly
compact set.

Assume now that $\phi_i(a)$ is compact for all $i\in I$. Since $\phi(a)$ is
not compact there exists an
$\varepsilon>0$ such that the set $\{i\in
I:\|\phi_i(a)\|\geq\varepsilon\}$ is infinite. Then the set $\sum_{i\in I}\oplus
\phi_i(a)\mathcal{B}(H_i)_{1/2}\phi_i(a)$ is not compact since it contains a copy of an
$l^{\infty}$
ball \cite[Chapter V, Theorem 4.2]{cnw}. This completes the proof as the last
set is
contained in
$\cp^2(a)$.
\end{proof}

\begin{remark}
 Let $\mathcal{A}$ be a C*-algebra and $a\in \mathcal{A}_1$. Then by the theorem above
and \cite[Theorem
2.2]{1996}, the following assertions are equivalent:
\begin{enumerate}
 \item There exists a faithful representation $(\phi,H)$ of $\mathcal{A}$ so that
$\phi(a)$ is a
compact operator.
 \item The set $\cp^2(a)$ is norm compact.
 \item The set $\cp^2(a)$ is weakly compact.
\end{enumerate}
\end{remark}

\vspace{1em}

 Let $\phi:\mathcal{V}\rightarrow \mathcal{B}(H_2,H_1)$ be a representation of a TRO
$\mathcal{V}$ and $K_1\subseteq H_1$ and
$K_2\subseteq H_2$ closed subspaces. A pair of subspaces $(K_2,K_1)$ is said to be
$\phi$-\textit{invariant} if $\phi(\mathcal{V})K_2\subseteq K_1$ and
$\phi(\mathcal{V})^*K_1\subseteq K_2$. The representation $\phi$ is said to be
\textit{irreducible} if $(0,0)$ and $(H_2,H_1)$ are the only
$\phi$-invariant pairs.

Two representations $\phi_i : \mathcal{V}\rightarrow \mathcal{B}(H_1^i,H_2^i)$ of
$\mathcal{V}$, $i=1,2$ are said to be unitarily equivalent, if there are unitary operators
$U_i: H_i^1\rightarrow H_i^2$, $i=1,2$ such that $\phi_1(x)=U_2^*\phi_2(x)U_1$, for all
$x\in\mathcal{V}$.

Let $(\phi_i)_{i\in I}$ be a maximal family of pairwise inequivalent irreducible
representations of $\mathcal{V}$, $\phi_i:\mathcal{V}\rightarrow
\mathcal{B}(H_{2,i},H_{1,i})$. Their direct sum
$\phi=\oplus_{i\in I} \phi_i$ is the reduced atomic representation of $\mathcal{V}$. It
follows from
\cite[Lemma 3.5]{ara} that an irreducible representation of a TRO is the
restriction of an irreducible representation of its linking algebra. Therefore, the
reduced atomic representation of a TRO $\mathcal{V}$ is the restriction of the reduced
atomic
representation of its linking algebra $A(\mathcal{V})$.

\begin{theorem}\label{gwctro} Let $\mathcal{V}$ be a TRO and $a\in \mathcal{V}_1$. The
following are
equivalent:
\begin{enumerate}
 \item $\cp^2_\mathcal{V}(a)$ is a weakly compact set.
 \item There exists a faithful representation $\pi$ of $\mathcal{V}$ such that $\pi(a)$
is a compact operator.
 \item $\phi(a)$ is a compact operator where $\phi$ is the reduced atomic representation
of $\mathcal{V}$.
\end{enumerate}
\end{theorem}
\begin{proof} First we show that (1) is equivalent to (2).
Suppose that the set $\cp^2_\mathcal{V}(a)$ is weakly compact. From Proposition
\ref{incl1} we know
that $\cp^2_{A(\mathcal{V})}(a)\subseteq \cp^2_\mathcal{V}(a)$ and
therefore the set $\cp^2_{A(\mathcal{V})}(a)$ is weakly compact. Now, by Theorem
\ref{cstarweakly}, there exists a faithful representation $\pi$ of
$\mathcal{V}$ that maps $a$ to a compact operator.

Conversely, suppose that $\pi$ is a faithful representation of $\mathcal{V}$ such that
$\pi(a)$ is
a compact
operator. We may assume that both $H_1$ and $H_2$ are infinte dimensional Hilbert spaces.
Identifying $\mathcal{V}$ with $\pi(\mathcal{V})$, Proposition \ref{incl}
states that 
$\cp_\mathcal{V}^2(a)\subseteq \mathcal{V}\cap \mathcal{K}(H_2,H_1)\subseteq
\mathcal{K}(H_1\oplus H_2)$. The set $\cp^2_\mathcal{V}(a)$ is
WOT-closed. Since the relative w* and WOT-topologies on the closed unit ball of
$\mathcal{B}(H_1\oplus H_2)$ coincide, \cite[Theorem 4.2.4.]{mur}, $\cp^2_\mathcal{V}(a)$
is a w*-closed
set. From the Banach-Alaoglu theorem we deduce that $\cp^2_\mathcal{V}(a)$ is a w*-compact
set.
By \cite[Proposition 10.4.3]{kr}, the weak topology on $\mathcal{K}(H_1\oplus H_2)$
coincides with
the relative w*-topology on $\mathcal{K}(H_1\oplus H_2)$ and therefore
$\cp^2_\mathcal{V}(a)$ is a weakly
compact set.

Obviously (3) implies (2) since $\phi$ is a faithful representation. So, we only need
to show that (1) implies (3). From Proposition \ref{incl1} we know that
$\cp_{A(\mathcal{V})}^2(a)\subseteq
\cp_\mathcal{V}^2(a)$. Therefore the set $\cp_{A(\mathcal{V})}^2(a)$ is weakly compact and
from
Theorem \ref{cstarweakly} $\rho(a)$ is a compact operator, where $\rho$ is the reduced
atomic representation of $A(\mathcal{V})$ \cite[Theorem 2.2]{1996}. The operator $\phi(a)$
is
compact since $\phi=\rho|_\mathcal{V}$.
\end{proof}

\vspace{1em}
Statement (3) of Theorem \ref{gwctro} ensures that the elements of $\mathcal{V}$ that are
mapped to
a compact operator by a faithful representation of $\mathcal{V}$ form a subTRO.

\begin{remark} \label{ccompact}
Let $\mathcal{A}$ be a C*-algebra which acts on a Hilbert space $H$ and contains
$\mathcal{K}(H)$, the set of
compact operators on $H$. If
$a\in \mathcal{A}_1$, the following assertions are equivalent:
\begin{enumerate}
 \item $a$ is a compact operator.
 \item The set $\cp^2(a)$ is norm compact.
 \item The set $\cp^2(a)$ is SOT-compact.
 \item The set $\cp^2(a)$ is weakly compact.
\end{enumerate}
\end{remark}
\begin{proof}
 From \cite[Corollary 2.4.]{1996}, we know that (1) and (2) are equivalent. 

That (2) implies (3) is obvious.

Now we show that (3) implies (1). The following arguments are similar to those
of \cite[Lemma 2.1]{1996}. Since
$\cp^2(a)$ is SOT-compact and $a(\mathcal{K}(H))_{1/2}a\subseteq\cp^2(a)$, the set
$a(\mathcal{K}(H))_{1/2}a$
is SOT-precompact. Let $\{f_n\}_{n=1}^{\infty}$ be
a bounded sequence in $H$. Without loss of generality we may assume that $\|f_n\|\leq
1/2$, for all $n\in\mathbb N$. Let $e$ be a unit
vector in $(\ker a^*)^{\perp}$. For every $n\in \mathbb N$, let $x_n=e\otimes
f_n$. Then
$ax_na=a^*e\otimes af_n$. Since $a(\mathcal{K}(H))_{1/2}a$ is
a SOT-precompact set, the sequence $\{(a^*e\otimes af_n)(h)\}_{n\in \mathbb N}$ has a
convergent subsequence for every $h\in H$. thus, the
sequence $\{\langle h, a^*e\rangle af_n\}_{n\in \mathbb N}$ has a convergent subsequence
and therefore $\{af_n\}_{n\in\mathbb N}$ has a
convergent subsequence. Hence, $a$ is a compact operator.

Obviously (1) implies (4). So, it suffices to see that (4) implies (1). Let us assume
that $a$ is a non-compact operator in $\mathcal{A}$. Following
the arguments of the proof of Proposition \ref{bofh} we can easily see that there exists
an
$\varepsilon>0$, an infinite dimensional
projection $p$ and a non-compact partial isometry $v$ on $p(H)$ such that
$\mathcal{K}(p(H),vp(H))_{\varepsilon^2/2}\subseteq \cp^2_{\mathcal{A}}(a)$. Since the
ball $\mathcal{K}(p(H),vp(H))_{\varepsilon^2/2}$ is not weakly compact, the set
$\cp^2_{\mathcal{A}}(a)$ is
not weakly compact.
\end{proof}

The following example shows that the compactness of
an element $u$ of a TRO does not imply the SOT-compactness of $\cp^2(u)$.

\begin{example} Let $\mathcal{V}$ be the TRO $\mathcal{B}(H_2,H_1)$ where $H_1$ is an
infinite dimensional
Hilbert space  and $H_2$ a one dimensional Hilbert space. The unit ball of
$\mathcal{B}(H_2,H_1)$
is not SOT-compact. Indeed, if $\{e\otimes
f_n\}_{n\in\mathbb N}$ is a sequence in $\mathcal{B}(H_2,H_1)$, where
$e$ is a unit vector of $H_2$ and $\{f_n\}$ an orthonormal sequence of $H_1$, then the
sequence $\{(e\otimes f_n)(e)\}=\{\langle e,e\rangle f_n\}=\{f_n\}$ has not a
convergent
subsequence. Consider an isometry $u\in \mathcal{V}$. Then $u$ is compact and $\cp^2(u)=\mathcal{V}_1$ is not a
SOT-compact set.
\end{example}

\vspace{1em}
The set $\cp^2(a)$ of the remark above is always WOT-compact since the WOT-topology of
$\mathcal{B}(H)$ coincides with its w*-topology on its closed unit ball. Therefore, the
WOT-compactness of $\cp^2(a)$ can not be equivalent with the statements of Remark
\ref{ccompact}.

\vspace{1em}
Vala introduced the notion of compactness in a normed algebra in \cite{1967}. He
defined an element $a$ of a normed algebra to be compact if the mapping $x \rightarrow
axa$ is compact.

\begin{definition}
 A linear mapping $u:\mathcal{V}\rightarrow \mathcal{V}$ is called a \textit{weakly
compact operator} on $\mathcal{V}$ if $\{u(x):\ \|x\|\leq 1\}$
is relatively weakly compact in $\mathcal{V}$.
\end{definition}

We shall use the following theorem. It was proved by K. Ylinen in \cite{1972} and
\cite{1975}.

\begin{theorem} \label{ylinenweak} 
Let $a$ be an element of the C*-algebra $\mathcal{A}$. The following conditions are
equivalent.
\begin{enumerate}
 \item There exists a faithful representation $\phi$ that maps $a$ to a compact operator.
 \item The operator $u:\mathcal{V}\rightarrow \mathcal{V}$, $u(x)=axa$ is compact.
 \item The operator $u:\mathcal{V}\rightarrow \mathcal{V}$, $u(x)=axa$ is weakly
compact.
\end{enumerate}
\end{theorem}

Bunce and Chu in \cite{1992} establish several theorems classifying
compact and weakly compact JB*-triples. A JB*-triple $\mathcal{A}$ is said to be (weakly)
compact if
the antilinear operator $x \rightarrow \{axa\}$ is (weakly) compact for each
$a\in\mathcal{A}$, where $\{\ \}$ denotes the ternary product. It
follows from \cite[Theorem 3.6]{1992} that a TRO $\mathcal{V}$ is isomorphic to a subTRO
of
$\mathcal{K}(H)$ for some Hilbert space $H$ if and only if the mapping $a\rightarrow
ax^*a$ is
compact or equivalently weakly compact, for all $a\in \mathcal{V}$. Next theorem
characterizes the compact
elements of a TRO $\mathcal{V}$.

\begin{theorem}\label{vala} Let $a$ be an element of a TRO $\mathcal{V}$. The following
conditions are equivalent:

\begin{enumerate}
 \item There exists a faithful representation $\pi$ that maps $a$ to a compact operator.
 \item The operator $u:\mathcal{V}\rightarrow \mathcal{V}$, $u(x)=ax^*a$ is compact.
 \item The operator $u:\mathcal{V}\rightarrow \mathcal{V}$, $u(x)=ax^*a$ is weakly
compact.
\end{enumerate}
\end{theorem}
\begin{proof} 
First we show that (2) implies (1). Let $u:\mathcal{V}\rightarrow \mathcal{V}$,
$u(x)=ax^*a$ be a compact
operator. Then the extension of $u$ to the linking algebra $A(\mathcal{V})$ of
$\mathcal{V}$ is compact as
well, since 
$\left(\begin{array}{cc}
    0 & a\\
    0 & 0
   \end{array} \right)
\left(\begin{array}{cc}
    x_1 & x_2\\
    x_3 & x_4
   \end{array} \right)
\left(\begin{array}{cc}
    0 & a\\
    0 & 0
   \end{array} \right)=
\left(\begin{array}{cc}
    0 & ax_3a\\
    0 & 0
   \end{array} \right)\in \mathcal{V}$,
where 
$\left(\begin{array}{cc}
    x_1 & x_2\\
    x_3 & x_4
   \end{array} \right) \in A(\mathcal{V}).$
It follows that the operator $\tilde{u}:A(\mathcal{V})\rightarrow A(\mathcal{V})$,
$\tilde{u}(x)=axa$ is compact. From \cite{1972} there exists a faithful
representation $\pi$ of $A(\mathcal{V})$ such that $\pi(a)$ is a compact operator.

Now we show the implication (1)$\Rightarrow$(2). Suppose there exists an isometric
representation $\pi$ of $\mathcal{V}$ on a Hilbert space $H$ so that $\pi(a)$ is a compact
operator
on
$H$. Then (see \cite{1964}) the map  $u_1:\mathcal{B}(H)\rightarrow \mathcal{B}(H)$,
$u_1(x)=\pi(a)x\pi(a)$ is compact.
Obviously, the map $u_2:\mathcal{B}(H)\rightarrow \mathcal{B}(H)$,
$u_2(x)=\pi(a)x^*\pi(a)$ is compact as
well. Therefore, the restriction of $u_2$ to $\pi(\mathcal{V})$ is a compact operator.
Since $\pi$ is an isometry the result follows.

That (1) implies (3) can be readily verified. 

Applying the arguments at the beginning of this proof and Theorem \ref{ylinenweak} we
deduce that (3) implies (1).
\end{proof}

\begin{remark}
 Let $\mathcal{V}$ be a TRO. It follows from Remark \ref{strict} and Theorem \ref{gwctro}
that the weak compactness of $\cp_\mathcal{V}^2(a)$ does not
imply its norm compactness. On the other hand, we would like to note that the norm
compactness and weak compactness of mapping $u:\mathcal{V}\rightarrow
\mathcal{V}$, $u(x)=ax^*a$ are equivalent.
\end{remark}

\vspace{1em}

\noindent \textbf{Acknowledgements.} I would like to thank Prof. M. Anoussis,
whose help, support and suggestions have been instrumental in completing this paper. I
would like also to thank the referee for his useful suggestions.

\bibliographystyle{amsplain}

\begin{thebibliography}{20}
\bibitem{2008} M. Anoussis, V. Felouzis and I. G. Todorov, \textit{Contractive perturbations in C*-algebras}, J.
Operator Theory \textbf{59}
(2008), no. 1, 53-68 MR2404464 (2009h:46100).
\bibitem{1996} M. Anoussis and E. G. Katsoulis, \textit{Compact operators and the
geometric structure of C*-algebras}, Proc. Amer. Math. Soc. \textbf{124} (1996),
2115-2122, MR1322911 (96i:46068).
\bibitem{1997} M. Anoussis and E. G. Katsoulis, \textit{Compact operators and the
geometric structure of nest algebras}, Indiana Univ. Math. J. \textbf{46} (1997), 319-335, MR1444482 (98e:47066a).
\bibitem{2005} M. Anoussis and I. Todorov, \textit{Compact operators on Hilbert modules},
Proc. Amer. Math. Soc. \textbf{133} (2005), 257-261, MR2086218 (2005e:46091).
\bibitem{ara} L. Aramba\v{s}i\'{c}, \textit{Irreducible representations of Hilbert
C*-modules}, Math. Proc. R. Ir. Acad. \textbf{105}A (2005), 11-24, MR2162903 (2006d:46078).
\bibitem{1992} L. J. Bunce and C.-H. Chu, \textit{Compact operations, multipliers and Radon-Nikodym property in
JB*-triples}, Pacific J. Math. \textbf{153} (1992), no. 2, 249-265, MR1151560 (93c:46125).
\bibitem{cnw} J.B. Conway, \textit{A Course in Functional Analysis}, Graduate Texts in Mathematics 96, Springer-Verlag,
New York,
(1990), MR1070713 (91e:46001).
\bibitem{1981} L. Harris, \textit{A generalization of C*-algebras}, Proc. London Math.
Soc. \textbf{42} (1981), 331-361, MR0607306 (82e:46089).
\bibitem{hes} M. Hestenes, \textit{A ternary algebra with applications to matrices and
linear transformations}, Arch.
Rational Mech.
Anal. \textbf{11} (1962), 138-194, MR0150166 (27 \#169).
\bibitem{kr} R. V. Kadison and J. R. Ringrose, \textit{Fundamentals of the theory of
operator algebras II}. Academic Press, New York, (1986), MR1468230 (98f:46001b).
\bibitem{2004} E. G. Katsoulis, \textit{Geometry of the unit ball and representation
theory for operator algebras}, Pacific J. Math. \textbf{216} (2004), 267-292, MR2094546
(2006b:47123).
\bibitem{1987} R. L. Moore and T. T. Trent, \textit{Extreme point of certain operator
algebras}, Indiana Univ. Math. J.
\textbf{36} (1987),
645-650, MR0905616 (89d:47103).
\bibitem{mur} G. J. Murphy, \textit{C*-algebras and operator theory}, Academic Press, (1990), MR1074574 (91m:46084).
\bibitem{1964} K. Vala, \textit{On compact sets of compact operators}, Ann. Acad. Sci.
Fenn. Ser. A I No. \textbf{351} (1964), MR0169078 (29 \#6333).
\bibitem{1967} K. Vala, \textit{Sur les \'{e}l\'{e}ments compacts d'une alg\`{e}bre
norm\'{e}e}, Ann. Acad. Sci. Fenn. Ser. A I No. \textbf{407} (1967), MR0222642 (36
\#5692).
\bibitem{1972} K. Ylinen, \textit{A note on the compact elements of C*-algebras}, Proc.
Amer. Math. Soc. \textbf{35} (1972), 305-306, MR0296716 (45 \#5775).
\bibitem{1975} K. Ylinen, \textit{Dual C*-algebras, weakly semi-completely continuous
elements, and the extreme rays of the positive cone}, Ann. Acad. Sci. Fenn. Ser. A I
Math. No \textbf{599} (1975), 9pp, MR0385584 (52 \#6445).
\bibitem{zettl} H. Zettl, \textit{A characterization of ternary rings of operators}, Adv. in Math. \textbf{48} (1983),
no. 2, 117-143, MR0700979 (84h:46093).
\end{thebibliography}

\end{document}